\documentclass[11pt]{article}

\usepackage{amsfonts,amsmath}
\usepackage{latexsym}

\author{K\'aroly J. B\"or\"oczky\footnote{Supported by
NKFIH grant K 132002}\footnote{Alfr\'ed R\'enyi Institute of Mathematics, 
 Realtanoda u. 13-15, H-1053 Budapest, Hungary, and
Department of Mathematics, Central European University, Nador u. 9, H-1051, Budapest, Hungary,
 boroczky.karoly.j@renyi.hu}, Apratim De\footnote{Department of Mathematics, Central European University, Nador u. 9, H-1051, Budapest, Hungary,
de.apratim91@gmail.com}}
\title{Stable solution of the Logarithmic Minkowski problem in the case of hyperplane symmetries}

\newcommand{\proof}{\noindent{\it Proof: }}
\newcommand{\proofbox}{\mbox{ $\Box$}\\}
\newcommand{\R}{\mathbb{R}}

\newtheorem{lemma}{LEMMA}[section]
\newtheorem{theo}[lemma]{THEOREM}
\newtheorem{example}[lemma]{EXAMPLE}

\newtheorem{coro}[lemma]{COROLLARY}
\newtheorem{conj}[lemma]{CONJECTURE}
\newtheorem{prop}[lemma]{PROPOSITION}

\begin{document}

\maketitle

\begin{abstract}
In the case of symmetries with respect to a
Coxeter group $G\subset O(n)$ acting without non-zero fixed points on $\R^n$, the stability
of the solution of the
  Logarithmic Minkowski problem on $S^{n-1}$ is established.
\end{abstract}


\section{Introduction}

The so called Minkowski problem form the core of various areas in fully nonlinear partial differential equations and convex geometry
(see Trudinger, Wang \cite{wang} and Schneider \cite{Sch14}), which was
extended to the $L_p$-Minkowski theory
 by Lutwak \cite{Lut93,Lut93a,Lut96}.
The classical Minkowski's existence theorem due to Minkowski and Aleksandrov describes the so called surface area measure
$S_K$ of a convex body $K$ (the case $p=1$)  where the regularity of the solution  is well investigated by
 Nirenberg \cite{NIR}, Cheng and Yau \cite{CY}, Pogorelov \cite{POG} and Caffarelli \cite{LC}.

First major results about the $L_p$-Minkowski problem for $p\neq 1$ have been obtained by
Chou, Wang \cite{CW} and
Hug, Lutwak, Yang, Zhang \cite{HLYZ2}, and more recently the papers
Boroczky,  Lutwak,  Yang,  Zhang \cite{BLYZ13},
Andrews, Guan, Ni \cite{AGN16},
Kolesnikov, Milman \cite{KolMilsupernew},
Bianchi, Boroczky, Colesanti, Yang \cite{BBCY}, Chen, Li, Zhu \cite{CLZ17},
Chen, Huang, Li \cite{CHL20}, Bryan, Ivaki, Scheuer \cite{BIJ19},
 present new developments and approaches. Cases of multiple solutions are discussed in
\cite{BLYZ13},  \cite{CLZ17}, He, Li, Wang \cite{HLW16} and
Li \cite{Li19}.

For a compact convex set $K$ in $\R^n$, we write $V(K)$ to denote to $n$-dimensional Lebesgue measure.
We say that a compact compact set  $K$ in $\R^n$ is a convex body if $V(K)>0$; or equivalently, the interior of $K$ is non-empty.
The cone volume measure or $L_0$-surface area measure $V_K$ on  $S^{n-1}$, 
whose study was initiated independently by Firey \cite{Fir74}
and Gromov and Milman \cite{GromovMilman},  has become an indispensable tool in the last decades (see say Barthe, Gu\'{e}don, Mendelson, Naor \cite{BG}, Naor \cite{Nar07},  Paouris, Werner \cite{PW}, Boroczky, Henk \cite{BoH16}).
If a convex body $K$ contains the origin, then its cone volume measure is
$dV_K=\frac1n\,h_K\,dS_K$ where $h_K$ is the support function of $K$ and the total measure is the volume of $K$.
In particular, the Monge-Amp\`ere equation on the sphere $S^{n-1}$ corresponding to the logarithmic (or $L_0$-) Minkowski problem is
\begin{equation}
\label{MongeVK}
h\det(\nabla^2 h+h\,{\rm Id})=  nf
\end{equation}
where $\nabla h$ and  $\nabla^2 h$ are the gradient and the Hessian of $h$ with respect to a moving orthonormal frame. We recall that for a given finite Borel measure $\mu$ on $S^{n-1}$, a positive $h$ on $S^{n-1}$ that is the restriction of a convex homogeneous function on $\R^n$ is the solution of \eqref{MongeVK} in the Alexandrov sense if the corresponding
Monge-Amp\`ere measure satisfies
\begin{equation}
\label{MongeVK-Alexandrov}
\det(\nabla^2 h+h\,{\rm Id})\,d\sigma=  \frac{n}{h}\cdot \mu
\end{equation}
where $\sigma$ is the Lebesgue measure on $S^{n-1}$.
In particular, for any Alexandrov solution $h$ of \eqref{MongeVK} (or equivalently of
\eqref{MongeVK-Alexandrov}), there exists a unique convex body $K$ with $o\in{\rm int}\,K$ such that
$h=h_K$ where $h_K(u)=\max_{x\in K}\langle u,x\rangle$ is the support function of $K$ for $u\in\R^n$,
$\mu=V_K$ is the cone-volume measure of $K$
and the corresponding Monge-Amp\`ere measure is $S_K$.

We observe that the Monge-Amp\`ere equation
\eqref{MongeVK} is homogeneous in the sense that replacing $f$ by $\lambda\,f$ for $\lambda>0$ is equivalent replacing $h$ by $\lambda^{1/n}h$. Therefore, we may assume that $V(K)=1$; or in other words, the $f$ in \eqref{MongeVK}
is a probability density, or the $\mu$ in \eqref{MongeVK-Alexandrov} is a probability measure.

Following partial and related results by Andrews \cite{And99}, Chou, Wang \cite{CW},
He, Leng, Li \cite{HLL06},
Henk, Sch\"urman, Wills \cite{HSW06}, Stancu \cite{Stancu},
Xiong \cite{Xio10}
the paper Boroczky, Lutwak, Yang, Zhang \cite{BLYZ13} characterized even cone volume measures
by the so called subspace concentration condition. Recently, breakthrough
results have been obtained by
Chen, Li, Zhu \cite{CLZ19}, Chen, Huang, Li \cite{CHL20}, Kolesnikov \cite{Kol20},
Nayar, Tkocz \cite{NaT20},
Kolesnikov, Milman \cite{KolMilsupernew},
 Putterman \cite{Put21} about the uniqueness of the solution, which is intimately related to the conjectured log-Minkowski inequality
 Conjecture~\ref{logMconj}. As it turns out, subspace concentration condition also holds for the cone-volume measure $V_K$ if
the centroid of a general convex body $K$ is the origin (see Henk, Linke \cite{HeL14} and B\"or\"oczky, Henk \cite{BoH16,BoH17}).

We note that the conjectured uniqueness of the solution  of the Logarithmic, or $L_0$-Minkowski problem \eqref{MongeVK} for even positive $C^\infty$ $f$ has a special role within the $L_p$-Minkowski Problems as if $p<0$, then
it is known that the solution may not be unique (see
Jian, Lu, Wang \cite{JLW15},  Li, Liu, Lu \cite{LLL}, Milman \cite{Mil}). On the positive side, extending the work
of Kolesnikov, Milman \cite{KolMilsupernew},
Chen, Huang, Li, Liu \cite{CHL20} verified the uniqueness of the solution if $1-\frac{c}{n^{3/2}}<p<1$
for some absolute constant $c\in(0,1)$ (see also Putterman \cite{Put21}).

Concerning possibly non-even measures, the  logarithmic Minkowski problem \eqref{MongeVK-Alexandrov} is wide open, as the best sufficient condition for a measure being a cone-volume measure is provided by Chen, Li, Zhu \cite{CLZ19} (solving for example the case of absolutely continuous measures), and some obstruction (necessary condition) is provided by
Boroczky, Hegedus \cite{BoH15}.

Boroczky, Kalantzopoulos \cite{BoK} proved the following characterization of cone-volume measures under hyperplane symmetry assumption. We note that for any group $G\subset O(n)$ acting on $\R^n$ without non-zero fixed points, there exist only finitely many $G$ invariant linear subspaces of $\R^n$ where $G$ is a Coxeter group if it is generated by reflections through $n$ independent hyperplanes.

\begin{theo}[Boroczky, Kalantzopoulos]
\label{VKsymchar}
Let $G\subset O(n)$ be a Coxeter group acting on $\R^n$ without non-zero fixed points.
For a finite non-trivial Borel measure $\mu$ on $S^{n-1}$ invariant under $G$,
there exists a $G$ invariant Alexandrov solution
of the logarithmic Minkowski equation \eqref{MongeVK-Alexandrov}
 if and only if
\begin{description}
\item[(i)] $\mu(L\cap S^{n-1})\leq \frac{{\rm dim}\,L}{n}\cdot\mu(S^{n-1})$ for any proper linear subspace $L$
invariant under $G$;
\item[(ii)] $\mu(L\cap S^{n-1})=\frac{{\rm dim}\,L}{n}\cdot\mu(S^{n-1})$ in (i)
for a proper invariant linear subspace $L$ is equivalent with
${\rm supp}\,\mu\subset L\cup L^\bot$.
\end{description}
In addition, if strict inequality holds in (i) for  each proper linear subspace $L$
invariant under $G$, then the $G$ invariant solution is unique.
\end{theo}

We note that the measure in Theorem~\ref{VKsymchar} may not be even; for example, possibly $\mu=V_K$ for a regular simplex $K$ whose centroid is the origin.

For compact convex sets $M$ and $N$, we write $M\oplus N$ to denote $M+N$
if $\langle x,y\rangle=0$ holds for $x\in M$ and $y\in N$. In addition, we say that a linear subspace $L$ of $\R^n$ is proper if
$1\leq {\rm dim}\, L\leq n-1$.
We note that \cite{BoK} proved that
$V_K(L\cap S^{n-1})=\frac{{\rm dim}\,L}{n}\cdot V(K)$
holds in Theorem~\ref{VKsymchar} (i) for a proper invariant subspace $L$ if and only if
$K=(K\cap L)\oplus(K\cap L^\bot)$.

According to \cite{BoK}, $V_K=V_C$ holds
for convex bodies $K$ and $C$  in $\R^n$ invariant under 
 a Coxeter group $G\subset O(n)$ acting on $\R^n$ without non-zero fixed points
 if and only if $V(K)=V(C)$, and $K=K_1\oplus\ldots \oplus K_m$ and $C=C_1\oplus\ldots \oplus C_m$ for compact convex sets $K_1,\ldots, K_m,C_1,\ldots,C_m$ of dimension at least one and invariant
under $G$
where $K_i$ and $C_i$ are dilates for $i=1,\ldots,m$.
Naturally, if $m=1$, then $K=C$.\\

In order to prepare for the stability version Theorem~\ref{VKVCcloseh} of Theorem~\ref{VKsymchar},
for any compact $X\subset S^{n-1}$ and $\varrho\in[0,2]$, we consider the tube
$$
\Psi(X,\varrho)=\{u\in S^{n-1}:\,\exists x\in X,\;\|x-u\| \leq \varrho\}.
$$

The cone volume measure $V_K$ of a convex body $K$ readily satisfies
$dV_{tK}=t^ndV_K$ for $t>0$. Therefore, when comparing the cone volume measures of convex bodies $K$ and $C$, we may asssume that $V(K)=V(C)=1$, and hence $V_K$ and $V_C$ are probability measures on $S^{n-1}$. In turn, one natural distance between two probability measures $\mu$ and $\nu$ on $S^{n-1}$ is the $l_1$ Wasserstein distance.
First, we consider the family of Lipschitz functions on $S^{n-1}$; namely, for $\theta>0$, let
\begin{equation}
\label{Wass-def}
{\rm Lip}_\theta=\big\{f:\,S^{n-1}\to\R:\,\forall a,b\in S^{n-1},\;|f(a)-f(b)|\leq\theta\|a-b\|\big\}.
\end{equation}
Now the Wasserstein distance of the Borel probability measures $\mu$ and $\nu$ on $S^{n-1}$ is
$$
d_W(\mu,\nu)=\sup\left\{\int_{S^{n-1}} f\,d\mu-\int_{S^{n-1}} f\,d\nu:\,f\in {\rm Lip}_1\right\}.
$$
It is known that convergence of a sequence of probability measures with respect to the Wasserstein distance is equivalent with weak convergence.

We note that as $\mu(S^{n-1})=\nu(S^{n-1})$ in the definition of $d_W(\mu,\nu)$, we may assume that
$\min f=-1$; therefore, $f\in {\rm Lip}_1$ implies that
\begin{equation}
\label{Linfinity-def}
\|f\|_\infty=\max_{u\in S^{n-1}}|f(u)|\leq 1.
\end{equation}
In turn, we observe that if $d\mu(u)=\varphi(u)\,du$ and $d\nu(u)=\psi(u)\,du$, then
\begin{equation}
\label{Wass-densfunction}
d_W(\mu,\nu)\leq \int_{S^{n-1}}|\varphi(u)-\psi(u)|\,du.
\end{equation}

\begin{theo}
\label{VKVCcloseh}
Let $G\subset O(n)$ be a Coxeter group acting on $\R^n$ without non-zero fixed points.
If $\mu_1$ and $\mu_2$ are Borel probability measures on $S^{n-1}$ invariant under
$G$, and
\begin{equation}
\label{VKVCclosesymcondh}
\begin{array}{rcl}
\mu_1\big(\Psi(L\cap S^{n-1},\delta)\big)&\leq& (1-\tau)\cdot\frac{{\rm dim}\,L}{n},\\[1ex]
\mu_2\big(\Psi(L\cap S^{n-1},\delta)\big)&\leq& (1-\tau)\cdot\frac{{\rm dim}\,L}{n}
\end{array}
\end{equation}
for $\delta,\tau\in(0,\frac12)$ and for any proper subspace $L$ invariant under $G$, then
the unique $G$ invariant Alexandrov solution $h_i$ of the logarithmic Minkowski problem
\eqref{MongeVK-Alexandrov}
for $\mu=\mu_i$, $i=1,2$, satisfies
\begin{eqnarray}
\label{VKVCclosehinfty}
\|h_1-h_2\|_\infty&\leq & \gamma_0\cdot d_W(\mu_1,\mu_2)^{\frac1{95n}}\\[1ex]
\label{VKVCclosehradii}
r_0\leq &h_1,h_2&\leq R_0
\end{eqnarray}
where  for some absolute constant $c>1$, we have
\begin{itemize}
\item $R_0=n$, $r_0=\frac1e$, $\gamma_0=c^n $
and the condition \eqref{VKVCclosesymcondh} is irrelevant
provided the action of $G$ is irreducible;
\item
$R_0=\left(\frac{n^6}{\delta}\right)^{\frac1{\tau}}$,
$r_0= \frac{n^{\frac{n}2}}{5^n}\left(\frac{\delta}{n^6}\right)^{\frac{n-1}{\tau}}$
and $\gamma_0=\frac{c^{n}}{\tau}\cdot \delta^{\frac{-3n}{\tau}}n^{\frac{12n}{\tau}}$
provided the action of $G$ is reducible.
\end{itemize}
\end{theo}

Actually, Theorem~\ref{VKVCcloseh} can be extended to the case when $\mu_1(S^{n-1})\neq \mu_2(S^{n-1})$
(see Corollary~\ref{VKVCclosehneq}). In this case, we need the bounded Lipschitz distance
$d_{\rm bL}(\mu,\nu)$ of two Borel measures
$\mu$ and $\nu$ on $S^{n-1}$ (see Dudley \cite{Dud02}); namely,
$$
d_{\rm bL}(\mu,\nu)=\sup\left\{\int_{S^{n-1}} f\,d\mu-\int_{S^{n-1}} f\,d\nu:\,f\in {\rm Lip}_1
\mbox{ and }\|f\|_\infty\leq 1\right\}.
$$
Using the test function constant $1$ shows that
\begin{equation}
\label{bLdifference}
|\mu(S^{n-1})-\nu(S^{n-1})|\leq d_{\rm bL}(\mu,\nu).
\end{equation}
We observe that if $\mu(S^{n-1})=\nu(S^{n-1})=1$, then $d_{\rm bL}(\mu,\nu)=d_{W}(\mu,\nu)$.
On the other hand, if $\lambda>0$ and $\mu$ is any finite non-trivial Borel measure on $S^{n-1}$, then
\begin{equation}
\label{bLdifference-lambda}
 d_{\rm bL}(\mu,\lambda\mu)\leq|\lambda-1|\cdot \mu(S^{n-1}).
\end{equation}

\begin{coro}
\label{VKVCclosehneq}
Let $G\subset O(n)$ be a Coxeter group acting on $\R^n$ without non-zero fixed points.
If $\mu_1$ and $\mu_2$ are finite Borel measures on $S^{n-1}$ invariant under
$G$ satisfying $ d_{\rm bL}(\mu_1,\mu_2)\leq M=\min\{\mu_1(S^{n-1}),\mu_2(S^{n-1})\}>0$  and
\begin{equation}
\label{VKVCclosesymcondhneq}
\begin{array}{rcl}
\mu_1\big(\Psi(L\cap S^{n-1},\delta)\big)&\leq& (1-\tau)\cdot\frac{{\rm dim}\,L}{n},\\[1ex]
\mu_2\big(\Psi(L\cap S^{n-1},\delta)\big)&\leq& (1-\tau)\cdot\frac{{\rm dim}\,L}{n}
\end{array}
\end{equation}
for $\delta,\tau\in(0,\frac12)$ and for any proper subspace $L$ invariant under $G$, then
the unique $G$ invariant Alexandrov solution $h_i$  of the logarithmic Minkowski problem
\eqref{MongeVK-Alexandrov}
for $\mu=\mu_i$, $i=1,2$, satisfies
\begin{eqnarray}
\label{VKVCclosehinftyneq}
\|h_1-h_2\|_\infty&\leq & \gamma_0M^{\frac1n}\cdot d_{\rm bL}(\mu_1,\mu_2)^{\frac1{95n}}\\[1ex]
\label{VKVCclosehradiineq}
r_0M^{\frac1n}\leq &h_1,h_2&\leq R_0M^{\frac1n}
\end{eqnarray}
where  for some absolute constant $c>1$, we have
\begin{itemize}
\item $R_0=2n$, $r_0=\frac1e$, $\gamma_0=c^n$
and the condition \eqref{VKVCclosesymcondhneq} is irrelevant
provided the action of $G$ is irreducible;
\item
$R_0=2\left(\frac{n^6}{\delta}\right)^{\frac1{\tau}}$,
$r_0= \frac{n^{\frac{n}2}}{5^n}\left(\frac{\delta}{n^6}\right)^{\frac{n-1}{\tau}}$
and $\gamma_0=\frac{c^{n}}{\tau}\cdot \delta^{\frac{-3n}{\tau}}n^{\frac{12n}{\tau}}$
provided the action of $G$ is reducible.
\end{itemize}
\end{coro}

Geometric inequalities under $n$ independent hyperplane symmetries were first considered by
Barthe, Fradelizi \cite{BaF13} and Barthe, Cordero-Erausquin \cite{BaC13}. These papers verified
the classical Mahler conjecture and Slicing conjecture, respectively, for these type of bodies.

We observe that the error term in Theorem~\ref{VKVCcloseh} in terms of $\varepsilon$ is not far from being optimal.
We provide an unconditional example; namely,
when $G$ is generated by the reflections through the coordinate hyperplanes.
Let $K$ be the unit cube $K=[-\frac12,\frac12]^n$, and the unconditional $C$ be obtained from $K$ by chopping off vertices of $K$ using simplices of
volume $\varepsilon$ and rescaling (to ensure $V(C)=1$). Then $d_W(V_K,V_C)<\gamma_1\cdot\varepsilon$, while
$(1-\gamma_2\varepsilon^{\frac1{n}})K\not\subset C$ for suitable $\gamma_1,\gamma_2>0$ depending on $n$.

The stable solution
Theorem~\ref{VKVCcloseh} of the logarithmic Minkowski
problem under hyperplane symmetry  does use the metric structure on $S^{n-1}$.
The next example shows that we can't expect an
"affine invariant" stability version of Theorem~\ref{VKVCcloseh} even if the cone volume measure is affine invariant in certain sense. \\

\begin{example}
If $e\in S^{n-1}$, and $K$ and $C$ are any convex bodies in $\R^n$ containing the origin in their interior
with $V(K)=V(C)=1$ and $V_K(e^\bot\cap S^{n-1})=V_C(e^\bot\cap S^{n-1})=0$, and
$\Phi_s$ is the diagonal transformation with $\Phi_s(e)=s^{-(n-1)}e$ and  $\Phi_s(x)=sx$
for $x\in e^\bot$, then both $V_{\Phi_s K}$ and $V_{\Phi_s C}$ tend weakly to $\mu_0$ as $s$ tends to infinity
where $\mu_0$ denotes the  probability measure on $S^{n-1}$ with $\mu_0(\{\pm e\})=\frac12$.
In particular, $V_{\Phi_s K}$ and $V_{\Phi_s C}$ are arbitrarily close if $s$ is large.
\end{example}

Next, we consider  two partial converses of Theorem~\ref{VKVCcloseh} to show that concerning Theorem~\ref{VKVCcloseh},
both the conditions involved and the conclusion are of the right kind. The first result does not require any symmetry assumption.

\begin{theo}
\label{KCcloseh}
Let $\mu_1$ and $\mu_2$ be finite Borel measures on $S^{n-1}$
such that there exists Alexandrov solution $h_i$  of the logarithmic Minkowski problem
\eqref{MongeVK-Alexandrov}
for $\mu=\mu_i$ and $i=1,2$.
If $h_1,h_2<R$ for $R>0$,  then
$$
d_{\rm bL}(\mu_1,\mu_2)\leq \gamma(R,n)\cdot \sqrt{\|h_1-h_2\|_\infty}
$$
where $\gamma(R,n)>0$ depends on $R$ and $n$.
\end{theo}

Secondly, we show that  if we have almost equality in Theorem~\ref{VKsymchar} (ii) for
measures $\mu_1$ and $\mu_2$ and
a proper linear subspace $L$
invariant under reflections through independent hyperplanes $H_1,\ldots,H_n$, then even if $\mu_1$ and $\mu_2$ are close,
it is possible that the  solutions $h_1$ and
  $h_2$ of \eqref{MongeVK-Alexandrov} are arbitrarily  far away.

\begin{theo}
\label{VKclosemaxh}
Let $G\subset O(n)$ be a group acting without non-zero fixed points on $\R^n$, let $R>\sqrt{n}$,
and
let $h$ be a positive $G$ invariant Alexandrov solution of  \eqref{MongeVK-Alexandrov}
for a probability measure $\mu$ on $S^{n-1}$ with $h<R$ such that
$$
\mu\big(\Psi(L\cap S^{n-1},\delta)\big)\geq (1-\varepsilon)\cdot\frac{{\rm dim}\,L}{n}
$$
for $\varepsilon\in(0,\frac{\varepsilon_0}{R^n})$, $\delta\in (0,\varepsilon]$ and
a proper subspace $L$ invariant under $G$ where $\varepsilon_0>0$ depends on $n$.
Then for any $t>1$, there exists a positive $G$ invariant Alexandrov solution $h_t$ of  \eqref{MongeVK-Alexandrov}
for a probability measure $\mu_t$ on $S^{n-1}$ such that
\begin{eqnarray*}
\|h-h_t\|_\infty&\geq&t\\
d_W(\mu,\mu_t)&\leq&\gamma(R,n)\varepsilon^{\frac1{10n}}
\end{eqnarray*}
where $\gamma(R,n)>0$ depends on $R$ and $n$.
\end{theo}

Concerning the set-up of the paper,
Section~\ref{secdiameter} proves the lower and upper bounds \eqref{VKVCclosehradii} on $h_i$ in Theorem~\ref{VKVCcloseh}.
Next Section~\ref{seclogMinkowskiconj} reviews the logarithmic Minkowski conjecture whose stability version Theorem~\ref{logMstab}  in the case of convex bodies with many hyperplane symmetries is essential in proving Theorem~\ref{VKVCcloseh} in Section~\ref{secVKVCcloseh}, leading also to Corollary~\ref{VKVCclosehneq}.
Finally, the two partial converses Theorem~\ref{KCcloseh} and Theorem~\ref{VKclosemaxh}  of Theorem~\ref{VKVCcloseh}
are proved in  Section~\ref{secconverse}.

\section{Bounding the diameter of $K$ in terms of $V_K$}
\label{secdiameter}

 First we point out a simple relation for balls contained in and containing a convex body.

\begin{lemma}
\label{rRvolume}
If $K$ is a convex body in $\R^n$ whose centroid is the origin, and $K\subset R\,B^n$ for
$R>0$, then $rB^n\subset K$ for some
$$
r\geq  \frac{n^{\frac{n}2}}{5^n}\cdot \frac{V(K)}{R^{n-1}}.
$$
\end{lemma}
\proof We set $r>0$ be maximal with the property $rB^n\subset K$
Since the origin is the centroid of $K$, we have $-K\subset nK$, and hence
$K$ is contained in a cylinder whose height is $(n+1)r\leq 2nr$ and base is an
$(n-1)$-ball of radius $R$. Therefore,
$$
V(K)\leq 2n\kappa_{n-1}R^{n-1}r.
$$
As $\Gamma(t+1)>(\frac{t}{e})^t\sqrt{2\pi t}$ for $t\geq 1$ (see Artin \cite{Art15})
and $\kappa_{n-1}<\frac{\sqrt{n+1}}{\sqrt{2\pi}}\cdot \kappa_n$, we have
$$
\kappa_{n-1}<\frac{\sqrt{n+1}}{\sqrt{2\pi}}\cdot \kappa_n=
\frac{\sqrt{n+1}}{\sqrt{2\pi}}\cdot \frac{\pi^{\frac{n}2}}{\Gamma(\frac{n}2+1)}<
\frac{\sqrt{n+1}}{\sqrt{2\pi}}\cdot \frac{(e\pi)^{\frac{n}2}}{n^{\frac{n}2}\sqrt{\pi n}}
<\frac{e^{\frac{n}2}\pi^{\frac{n-2}2}}{n^{\frac{n}2}}.
$$
In turn, we deduce
$$
r\geq \frac{V(K)}{2n\kappa_{n-1}R^{n-1}}>
\frac{n^{\frac{n}2}}{n\cdot e^{\frac{n}2}\pi^{\frac{n}2}}\cdot \frac{V(K)}{R^{n-1}},
$$
completing the proof of Lemma~\ref{rRvolume} as $n\cdot e^{\frac{n}2}\pi^{\frac{n}2}<5^n$.
\proofbox

For a convex body $K$ in $\R^n$, we write $R(K)$ to denote the minimal radius  of a Euclidean ball containing $K$,
and $r(K)$ to denote the radius of largest ball contained in $K$.
We observe that if the convex body $K$ is invariant under the reflections through
the hyperplanes $H_1,\ldots,H_n$ with $H_1\cap \ldots\cap H_n=\{o\}$,
then its centroid is the origin, and
$$
r(K)B^n\subset K\subset R(K)\,B^n.
$$

For Proposition~\ref{logentropydiam} and Lemma~\ref{logentropydiamirred}, let $\widetilde{B}$ denote the Euclidean ball centered at the origin with $V(\widetilde{B})=1$.

\begin{prop}
\label{logentropydiam}
Let $n\geq 2$ and $\delta,\tau\in (0,\frac1{2})$,
let $G\subset O(n)$ be a Coxeter group acting reducibly and without non-zero fixed points on $\R^n$,
 and
let the Borel probability measure $\mu$ on $S^{n-1}$
be invariant under $G$ and satisfy
$$
\mu(\Psi(L\cap S^{n-1},\delta))<(1-\tau)\cdot \frac{i}{n}
$$
for any linear $i$-subspace $L$ of $\R^n$, $i=1,\ldots,n-1$,
invariant under $G$,
 and let $V(C)=1$ hold for convex body $C$ in $\R^n$ invariant under $G$. Then
\begin{description}
\item[(i)]
$$
\int_{S^{n-1}}\log h_C\,d\mu\geq \log\frac{R(C)^\tau\delta}{n^5}, \mbox{ \ }and
$$
\item[(ii)] if
$\int_{S^{n-1}}\log h_C\,d\mu\leq \int_{S^{n-1}}\log h_{\widetilde{B}}\,d\mu$, then
$$
R(C)<\left(\frac{n^6}{\delta}\right)^{\frac1{\tau}}\mbox{ \ and \ }
r(C)> \frac{n^{\frac{n}2}}{5^n}\left(\frac{\delta}{n^6}\right)^{\frac{n-1}{\tau}}.
$$
\end{description}
\end{prop}
\proof  Let $E$ be the John ellipsoid of maximal volume in $C$, and hence $E$ is invariant under $G$, and
\begin{equation}
\label{EK}
E\subset C\subset n\, E.
\end{equation}

Let $L_1,\ldots,L_m$ be the irreducible linear subspaces invariant under $G$. The symmetries
of $E$ yield that there exists a set of principal directions of $E$ that are part of $L_1\cup\ldots\cup L_m$, and
for each $L_i$ there exists $r_i>0$ such that $E\cap L_i=r_i(B^n\cap L_i)$, $i=1,\ldots,m$. We may assume that
$r_1\leq\ldots\leq r_m$.

If $m=1$, then \eqref{EK} yields that $r_1B^n\subset C\subset nr_1B^n$; therefore, Proposition~\ref{logentropydiam}
trivially holds. In particular, let
$$
m\geq 2.
$$

For
$$
Q={\rm conv}\{r_i(B^n\cap L_i)\}_{i=1,\ldots,m},
$$
$E$ is the so-called Loewner (minimal volume ellipsoid) of $Q$,
and hence $Q\subset E\subset \sqrt{n}\,Q$, thus
(\ref{EK}) yields that $Q\subset C\subset n^2 Q$. In particular,
writing $d_i={\rm dim}\,L_i$ for  $i=1,\ldots,m$,
$Q\subset C$ satisfies
\begin{equation}
\label{QK}
n^n\prod_{i=1}^m r_i^{d_i}\geq \prod_{i=1}^m r_i^{d_i}\kappa_{d_i}\geq V(Q)\geq n^{-2n}V(C)=n^{-2n}
\end{equation}
where $d_1+\ldots+d_m=n$.
We observe that for any $u\in S^{n-1}$ , there exists $L_i$ such that
$\|u|L_i\|\geq \frac{1}{\sqrt{m}}> \frac{\delta}{n}$.
For $i=1,\ldots,m$, we define
\begin{eqnarray*}
\Lambda_i&=&L_1\oplus\ldots\oplus L_i\\
B_i&=&\left\{u\in S^{n-1}:\|u|L_i\|\geq \frac{\delta}{n}\mbox{ and }
\|u|L_j\|<\frac{\delta}{n}\mbox{ for }j>i\right\}.
\end{eqnarray*}

It follows that $S^{n-1}$ is partitioned into the Borel sets $B_1,\ldots,B_m$, and
as $B_j\subset \Psi(\Lambda_i\cap S^{n-1},\delta)$ for $1\leq j\leq i\leq m-1$, we have
\begin{eqnarray}
\label{muAi}
\mu(B_1)+\ldots+\mu(B_i)&\leq &\frac{(d_1+\ldots+d_i)(1-\tau)}{n}\mbox{ \ for $i=1,\ldots,m-1$}\\
\label{muAn}
\mu(B_1)+\ldots+\mu(B_{m})&=&1.
\end{eqnarray}
For $\zeta=\frac{1-\tau}{n}>\frac1{2n}$, next we define
\begin{eqnarray}
\label{betai}
\beta_j&= &\mu(B_j)-d_j\zeta\mbox{ \ for $j=1,\ldots,m-1$}\\
\label{betan}
\beta_m&=&\mu(B_{m})-d_m\zeta-\tau
\end{eqnarray}
where (\ref{muAi}) and (\ref{muAn}) yield
\begin{eqnarray}
\label{sumbetai}
\beta_1+\ldots+\beta_i&\leq &0\mbox{ \ for $i=1,\ldots,m-1$}\\
\label{sumbetan}
\beta_1+\ldots+\beta_m&=&0.
\end{eqnarray}

It follows from $r_i\,B^n\cap L_i\subset Q$ and from the definition of $B_i$
 that
$h_Q(u)\geq r_i\cdot\frac{\delta}{n}$
for $u\in B_i$, $i=1,\ldots,m$.
We deduce from applying (\ref{QK}), (\ref{muAn}),  (\ref{betai}), (\ref{betan}),
(\ref{sumbetai}), (\ref{sumbetan}), $r_1\leq\ldots\leq r_{m}$
and $\frac1{2n}<\zeta<\frac1{n}$
that
\begin{eqnarray*}
\int_{S^{n-1}}\log h_C\,d\mu&\geq  &\int_{S^{n-1}}\log h_Q\,d\mu
=\sum_{i=1}^{m}\int_{B_i}\log h_Q\,d\mu\\
&\geq &\sum_{i=1}^{m}\mu(B_i)\log r_i+\sum_{i=1}^{m}\mu(B_i)\log \frac{\delta}{n}
=  \sum_{i=1}^{m}\mu(B_i)\log r_i+\log \frac{\delta}{n}\\
&=&
\sum_{i=1}^{m}\beta_i\log r_i+\sum_{i=1}^{m}\zeta d_i\log r_i
+\tau\log r_{m}+\log \frac{\delta}{n}\\
&\geq &
\sum_{i=1}^{m}\beta_i\log r_i+\zeta\log \frac{1}{n^{3n}}
+\tau\log r_{m}+\log \frac{\delta}{n}\\
&= &
(\beta_1+\ldots+\beta_{m})\log h_{m}+
\sum_{i=1}^{m-1}(\beta_1+\ldots+\beta_i)(\log r_i-\log r_{i+1})\\
&& -3n\zeta\log n
+\tau\log r_{m}+\log \frac{\delta}{n}\\
&\geq &-3\log n+\tau\log h_{n}+\log \frac{\delta}{n}
\end{eqnarray*}
where we used $\zeta<\frac1n$ at the end. Now
$r_m=R(E)\geq R(C)/n$ and $\tau<1$ imply
\begin{eqnarray*}
-3\log n+\tau \log r_{m}+\log \frac{\delta}{n}&\geq & -3\log n+\tau \log \frac{R(C)}{n}+\log \frac{\delta}{n}\\
&\geq& -3\log n+
\tau \log R(C) - \log n +\log\delta-\log n\\
&=&\log\frac{R(C)^\tau\delta}{n^5},
\end{eqnarray*}
proving Proposition~\ref{logentropydiam} (i).

For (ii), let $\tilde{r}_n$ be the radius of $\widetilde{B}$, and hence
$\Gamma(\frac{n}2+1)<(\frac{n}{2e})^{\frac{n}2}\sqrt{2\pi(\frac{n}2+1)}<(\frac{2n}{e})^{\frac{n}2}$ implies
$$
1=\tilde{r}_n^n \kappa_n=\tilde{r}_n^n\cdot\frac{\pi^{\frac{n}2}}{\Gamma(\frac{n}2+1)}>\tilde{r}_n^n\cdot
\left( \frac{e\pi}{2n}\right)^{\frac{n}2},
$$
and hence
\begin{equation}
\label{tildeBradius}
\tilde{r}_n<\sqrt{\frac{2n}{e\pi}}.
\end{equation}
We deduce from (i) and \eqref{tildeBradius} that
$$
\log\frac{R(C)^\tau\delta}{n^5}\leq \int_{S^{n-1}}\log h_C\,d\mu\leq
 \int_{S^{n-1}}\log h_{\widetilde{B}}\,d\mu<\log \sqrt{\frac{2n}{e\pi}},
$$
thus $R(C)<(\frac{n^6}{\delta})^{\frac1{\tau}}$.

In turn, the bound for $r(C)$ follows from Lemma~\ref{rRvolume},
completing the proof of Proposition~\ref{logentropydiam}.
\proofbox

\begin{lemma}
\label{logentropydiamirred}
Let the action of the Coxeter group $G\subset O(n)$ be irreducible,
 and
let the Borel probability measure $\mu$ on $S^{n-1}$
be invariant under $G$,
 and let $V(C)=1$ hold for convex body $C$ in $\R^n$ invariant under $G$.
Then
$$
\int_{S^{n-1}}\log h_C\,d\mu\geq -1;
$$
$$
\frac1{e}<r(C)\leq R(C)<n.
$$
\end{lemma}
\proof As the action of $G$ is irreducible, it follows that
the inscribed ball of $C$ is the John ellipsoid; namely, the ellipsoid of maximum volume contained in $C$. According to
Ball \cite{Ball91}, $r(C)$ is at least the inradius $r_n$ of the regular simplex of volume one,
and hence $n!\leq (\frac{n}e)^n\sqrt{2\pi\,n}$ (see Artin \cite{Art15})  yields
$$
r(C)^n\geq r_n^n=\frac{n!}{n^{\frac{n}2}(n+1)^{\frac{n+1}2}}\geq
\frac{(\frac{n}e)^n\sqrt{2\pi\,n}}{2n^{n+\frac12}}>\frac1{e^n}.
$$
On the other hand, as the action of $G$ is irreducible, it follows that
the circumscribed ball of $C$ is the Loewner ellipsoid; namely, the ellipsoid of minimum volume containing $C$. According to
Barthe \cite{Bar97} (see also Lutwak, Yang, Zhang \cite{LYZ07}), $R(C)$ is at most the inradius $R_n$ of the regular simplex of volume one,
and hence $n!< (\frac{n}e)^n\sqrt{2\pi(n+1)}$ (see Artin \cite{Art15})  yields
$$
R(C)^n\leq R_n^n=\frac{n^n\cdot n!}{n^{\frac{n}2}(n+1)^{\frac{n+1}2}} <
\frac{n^{\frac{n}2}\cdot (\frac{n}e)^n\sqrt{2\pi(n+1)}}{(n+1)^{\frac{n+1}2}}<
\frac{n^n\sqrt{2\pi}}{e^n}<n^n.
$$
We conclude $\frac1{e}<r(C)\leq R(C)<n$.

Finally, $r(C)>\frac1{e}$ implies that $\log h_C(u)>-1$ for all $u\in S^{n-1}$.
\proofbox

For a convex body $K$ with $V(K)=1$ and hyperplane symmetries, combining Proposition~\ref{logentropydiam} with the consequence
$\int_{S^{n-1}}\log h_K\,dV_K\leq \int_{S^{n-1}}\log h_{\widetilde{B}}\,dV_K$ of
the  Logarithmic Minkowski Inequality  Theorem~\ref{logMsym} or using Lemma~\ref{logentropydiamirred}
 yield the following.

\begin{coro}
\label{cone-volume-diam}
Let $G\subset O(n)$ be a Coxeter group acting without non-zero fixed points on $\R^n$, and
let the  convex body $K$ in $\R^n$, $n\geq 2$, be invariant under $G$.
If, for $\delta,\tau\in (0,\frac1{2})$ and $i=\{1,\ldots,n-1\}$,
$$
V_K(\Psi(L\cap S^{n-1},\delta))<(1-\tau)\cdot \frac{i}{n}\cdot V(K)
$$
for any $i$-dimensional   subspace $L$ invariant under $G$, then
\begin{eqnarray*}
R(K)&<& \left\{
\begin{array}{ll}
\left(\frac{n^6}{\delta}\right)^{\frac1{\tau}}V(K)^{\frac1n}&\mbox{ if the action of  $G$ is reducible;}\\[1ex]
nV(K)^{\frac1n}& \mbox{ if the action of  $G$ is irreducible;}
\end{array} \right.\\[1ex]
r(K)&>& \left\{
\begin{array}{ll}
\frac{n^{\frac{n}2}}{5^n}\left(\frac{\delta}{n^6}\right)^{\frac{n-1}{\tau}}
V(K)^{\frac1n}&\mbox{ if the action of  $G$ is reducible;}\\[1ex]
\frac1e\cdot V(K)^{\frac1n}& \mbox{ if the action of  $G$ is irreducible.}
\end{array} \right.
\end{eqnarray*}
\end{coro}

Another consequence of Proposition~\ref{logentropydiam} is a condition yielding that a convex body
with hyperplane symmetries
is not close to be the direct sum of lower dimensional invariant compact convex sets.

\begin{prop}
\label{noMinkowskisum}
Let $G\subset O(n)$ be a Coxeter group acting reducibly and without non-zero fixed points on $\R^n$,
and
let the  convex body $K$ in $\R^n$, $n\geq 2$, be invariant under
$G$.
If $\delta,\tau\in (0,\frac1{2})$,
and a convex body $K$ in $\R^n$ invariant under
 $G$ satisfies
$$
V_K(\Psi(L\cap S^{n-1},\delta))<(1-\tau)\cdot \frac{{\rm dim}L}{n}\cdot V(K)
$$
for any proper coordinate subspace $L$ invariant under $G$, then
$$
(1-\eta)\big((L\cap K)\oplus (L^\bot\cap K)\big)\not\subset K
$$
for any proper subspace $L$ invariant under $G$ where
$$
\eta = \frac{\delta\tau}{4n} \cdot
\frac{n^{\frac{n}2}}{5^n}\left(\frac{\delta}{n^6}\right)^{\frac{n}{\tau}}.
$$
\end{prop}
\proof We may assume that $V(K)=1$, and define
\begin{eqnarray*}
R_0&=& \left(\frac{n^6}{\delta}\right)^{\frac1{\tau}}\\
r_0&=& \frac{n^{\frac{n}2}}{5^n}\left(\frac{\delta}{n^6}\right)^{\frac{n-1}{\tau}},
\end{eqnarray*}
and hence
\begin{equation}
\label{etaR0r0fraction}
\eta= \frac{\delta\tau}{4n} \cdot\frac{r_0}{R_0}<\frac{\tau}{4n},
\end{equation}
while Proposition~\ref{logentropydiam} implies that
$$
r_0B^n\subset K\subset R_0 B^n.
$$

We prove Lemma~\ref{noMinkowskisum} by contradiction; therefore, we suppose that there exists
a coordinate $i$-subspace $L$, $1\leq i\leq n-1$, such that
\begin{equation}
\label{Lindirect}
(1-\eta)\big((L\cap K)\oplus (L^\bot\cap K)\big)\subset K.
\end{equation}
We define
$$
\Omega_0=\left\{[o,x+y]:\, x\in(1-\eta)\partial (L\cap K)\mbox{ \ and \ }
y\in (1-\eta)\left(1-\frac{\tau}{2n} \right) (L^\bot\cap K) \right\}.
$$
In addition, let
\begin{eqnarray*}
\Omega&=&\{z\in K:\, \exists t\in(0,1],\;tz\in \Omega_0\}\\
&=&\{z\in K:\, (1-\eta)x\in \Omega_0\}\\
\Xi&=&\{u\in S^{n-1}:\, \exists x\in \Omega\cap\partial K,\; h_C(u)=\langle x,u\rangle\}.
\end{eqnarray*}
We deduce using $\eta<\frac{\tau}{2n}$ that
\begin{eqnarray}
\nonumber
V_K(\Xi)&\geq &{\cal H}^n(\Omega_0)\\
\nonumber
&=&\frac{i}{n}\cdot(1-\eta)^i{\cal H}^i(L\cap K)
\cdot (1-\eta)^{n-i}\left(1-\frac{\tau}{2n} \right)^{n-i} {\cal H}^{n-i}(L^\bot\cap K)\\
&>&(1-\tau)\frac{i}{n}\cdot{\cal H}^i(L\cap K)\cdot {\cal H}^{n-i}(L^\bot\cap K)>
(1-\tau)\frac{i}{n}.
\end{eqnarray}
Therefore, we contradict \eqref{Lindirect} by proving
\begin{equation}
\label{XiinPsi}
\Xi\subset \Psi(L\cap S^{n-1},\delta).
\end{equation}
Let $u\in\Xi$ be an exterior normal at $z\in\partial K$.
We observe that
$$
u=v\cos\beta +w\sin \beta
$$
where $v\in L\cap S^{n-1}$, $w\in L\cup S^{n-1}$
and $\beta=\angle (u,v)\in[0,\frac{\pi}2)$.
We write
$z=x+y$ for $x\in L\cap K$ and $y\in L^\bot\cap K$.
As $z\in\Xi$, we have
\begin{eqnarray*}
(1-\eta)x+(1-\eta)y
&=&(1-\eta)z\in\Omega_0\\
&\in&(1-\eta)(L\cap K)+
(1-\eta)\left(1-\frac{\tau}{2n} \right) (L^\bot\cap K).
\end{eqnarray*}
In turn, we deduce that
\begin{equation}
\label{yinMinkowskisum}
y\in \left(1-\frac{\tau}{2n} \right) (L^\bot\cap K).
\end{equation}
Let
$$
p=(1-\eta)x+y+\frac{\tau}{4n}\cdot r_0 w,
$$
which, using \eqref{yinMinkowskisum}, $r_0B^n\subset K$, \eqref{etaR0r0fraction}
and \eqref{Lindirect} satisfies
\begin{eqnarray*}
p&\in&(1-\eta)(L\cap K)+\left(1-\frac{\tau}{2n} \right) (L^\bot\cap K)+
\frac{\tau}{4n}\cdot (L^\bot\cap K)\\
&=& (1-\eta)(L\cap K)+\left(1-\frac{\tau}{4n} \right) (L^\bot\cap K)\\
&\subset& (1-\eta)(L\cap K)+(1-\eta) (L^\bot\cap K)\subset K.
\end{eqnarray*}
Since $u$ is exterior normal at $z=x+y$ where
$w\in L^\bot\cap S^{n-1}$, $v\in L\cap S^{n-1}$ and $x\in L\cap R_0B^n$, we have
\begin{eqnarray*}
0&\geq &\langle u, p-z\rangle=\left\langle u, \frac{\tau r_0}{4n}\cdot  w-\eta x\right\rangle\\
&=&
\left\langle v\cos\beta +w\sin \beta, \frac{\tau r_0}{4n}\cdot  w-\eta x\right\rangle=
 \frac{\tau r_0}{4n}\cdot \sin\beta-\langle v,x\rangle\eta\cos\beta\\
&\geq& \frac{\tau r_0}{4n}\cdot \sin\beta-R_0\eta\cos\beta.
\end{eqnarray*}
We conclude that
$$
\|u-v\|\leq \tan \beta\leq \frac{4n\eta}{\tau}\cdot\frac{R_0}{r_0}\leq \delta,
$$
which in turn, yields \eqref{XiinPsi} and contradicts \eqref{Lindirect},
proving
Proposition~\ref{noMinkowskisum}.
\proofbox

\section{On the Logarithmic Minkowski conjecture}
\label{seclogMinkowskiconj}

For origin symmetric convex bodies, the following is an equivalent form of the  origin symmetric case of the Logarithmic Brunn-Minkowski conjecture (see Boroczky, Lutwak, Yang, Zhang \cite{BLYZ12}).

\begin{conj}[Logarithmic Minkowski conjecture]
\label{logMconj}
If $K$ and $C$ are convex bodies in $\R^n$ whose centroid is the origin, then
\begin{equation}
\label{logMconjeq}
\int_{S^{n-1}}\log \frac{h_C}{h_K}\,dV_K\geq \frac{V(K)}n\log\frac{V(C)}{V(K)}
\end{equation}
with equality if and only if $K=K_1+\ldots + K_m$ and $L=L_1+\ldots + L_m$  compact convex sets
	$K_1,\ldots, K_m,L_1,\ldots,L_m$ of dimension at least one where $\sum_{i=1}^m{\rm dim}\,K_i=n$
	and $K_i$ and $L_i$ are dilates, $i=1,\ldots,m$.
\end{conj}

The argument in Boroczky, Lutwak, Yang, Zhang \cite{BLYZ13} yields that 
uniqueness of the solution of the logarithmic-Minkowski problem \eqref{MongeVK} for any positive even $C^\infty$ $f$
is equivalent saying that the Logarithmic Minkowski conjecture \eqref{logMconjeq} holds
for any $o$-symmetric convex bodies $K$ and $C$ with $C^\infty_+$ boundaries with equality if and only if $K$ and $C$
are dilates.

In $\R^2$, Conjecture~\ref{logMconj} is verified in Boroczky, Lutwak, Yang, Zhang \cite{BLYZ12} for origin symmetric convex bodies, but it is still open in general.
In higher dimensions, Conjecture~\ref{logMconj} is proved
for with enough hyperplane symmetries ({\it cf.} Theorem~\ref{logMsym}) and
 complex bodies ({\it cf.} Rotem \cite{Rotem}).

For origin symmetric convex bodies, Conjecture~\ref{logMconj} is proved when $K$ is close to be an ellipsoid  by a combination of the local estimates by Kolesnikov, Milman \cite{KolMilsupernew} and the use of the continuity method in PDE by Chen, Huang, Li, Liu \cite{CHL20}. Another even more recent proof of this result based on Alexandrov's approach of considering the Hilbert-Brunn-Minkowski operator for polytopes  is due to Putterman \cite{Put21}.
Additional local versions of Conjecture~\ref{logMconj} are due to Kolesnikov, Livshyts \cite{KoL} and Hosle, Kolesnikov, Livshyts \cite{HKL}.

Following the result on unconditional convex bodies by Saroglou \cite{Sar15}, Boroczky, Kalantzopoulos \cite{BoK} verified the
logarithmic Minkowski conjecture for convex bodies with $n$ independent
  hyperplane symmetries.

\begin{theo}[Boroczky, Kalantzopoulos]
\label{logMsym}
If the convex bodies $K$ and $C$  in $\R^n$ are invariant under linear reflections $A_1,\ldots,A_n$ through $n$ independent
linear $(n-1)$-planes $H_1,\ldots,H_n$, then
$$
\int_{S^{n-1}}\log \frac{h_C}{h_K}\,dV_K\geq  \frac{V(K)}n\log\frac{V(C)}{V(K)},
$$
with equality
 if and only if $K=K_1+\ldots+ K_m$ and $L=L_1+\ldots + L_m$ for compact convex sets
	$K_1,\ldots, K_m,L_1,\ldots,L_m$ of dimension at least one and invariant under  $A_1,\ldots,A_n$
where $K_i$ and $L_i$ are dilates, $i=1,\ldots,m$, and $\sum_{i=1}^m{\rm dim}\,K_i=n$.
\end{theo}

The Boroczky, De \cite{BoD-BM} proved the following
 stability version of
the logarithmic-Minkowski inequality Theorem~\ref{logMsym}
 for convex bodies with many hyperplane symmetries.

\begin{theo}
\label{logMstab}
If  the convex bodies $K$ and $C$  in $\R^n$ are invariant under the Coxeter group
$G\subset O(n)$ acting without non-zero fixed points on $\R^n$, and
$$
\int_{S^{n-1}}\log \frac{h_C}{h_K}\,\frac{dV_K}{V(K)}\leq \frac{1}n\cdot \log\frac{V(C)}{V(K)}+\varepsilon
$$
for $\varepsilon>0$,  then  for some $m\geq 1$,
there exist  compact convex sets
$K_1,C_1,\ldots,K_m,C_m$  of dimension at least one and invariant under  $G$
where $K_i$ and $C_i$ are dilates, $i=1,\ldots,m$, and $\sum_{i=1}^m{\rm dim}\,K_i=n$ such that
\begin{eqnarray*}
K_1+\ldots + K_m\subset &K&\subset
\left(1+c^n\varepsilon^{\frac1{95n}}\right)(K_1+\ldots +K_m)\\
C_1+\ldots +C_m\subset &C&\subset
\left(1+c^n\varepsilon^{\frac1{95n}}\right)(C_1+\ldots +C_m)
\end{eqnarray*}
where $c>1$ is an absolute constant.
\end{theo}

\section{Proof of Theorem~\ref{VKVCcloseh}}
\label{secVKVCcloseh}

For compact convex sets $K$ and $C$ in $\R^n$, their Hausdorff distance is
$$
d_\infty(K,C)=\|h_K-h_C\|_\infty=\min\{r\geq 0:\,K\subset C+r\,B^n
\mbox{ and }C\subset K+r\,B^n\}.
$$
We prove Theorem~\ref{VKVCcloseh} in the following form.

\begin{theo}
\label{VKVCclose}
Let $G\subset O(n)$ be a Coxeter group acting without non-zero fixed points on $\R^n$.
If $K$ and $C$ are convex bodies in $\R^n$ invariant under $G$ and satisfy $V(K)=V(C)=1$,
\begin{equation}
\label{VKVCclosesymcond}
\begin{array}{rcl}
V_K\big(\Psi(L\cap S^{n-1},\delta)\big)&\leq& (1-\tau)\cdot\frac{{\rm dim}\,L}{n},\\[1ex]
V_C\big(\Psi(L\cap S^{n-1},\delta)\big)&\leq& (1-\tau)\cdot\frac{{\rm dim}\,L}{n}
\end{array}
\end{equation}
for $\delta,\tau\in(0,\frac12)$ and for any proper subspace $L$ invariant under $G$, then
\begin{eqnarray}
\label{r0hKCR0}
r_0<&h_K,\;h_C&<R_0;\\[1ex]
\label{KCdilate}
d_\infty(K,C)
&\leq &
\gamma_0\cdot d_W(V_K,V_C)^{\frac1{95n}}
\end{eqnarray}
where  for some absolute constant $c>1$, we have
\begin{itemize}
\item $R_0=n$, $r_0=\frac1e$ and $\gamma_0=c^n$ provided the action of $G$ is irreducible
(and hence the condition \eqref{VKVCclosesymcond} is irrelevant);
\item
$R_0=\left(\frac{n^6}{\delta}\right)^{\frac1{\tau}}$,
$r_0= \frac{n^{\frac{n}2}}{5^n}\left(\frac{\delta}{n^6}\right)^{\frac{n-1}{\tau}}$
and $\gamma_0=\frac{c^{n}}{\tau}\cdot \delta^{\frac{-3n}{\tau}}n^{\frac{12n}{\tau}}$
provided the action of $G$ is reducible.
\end{itemize}
\end{theo}

We need the simple statements Lemma~\ref{RKhKLipsitz} and \eqref{WassersteinLipsitz}.

\begin{lemma}
\label{RKhKLipsitz}
If $K$ is a convex body with $K\subset R\,B^n$ for $R>0$, then
$|h_K(u)-h_K(v)|\leq R\|u-v\|$ for $u,v\in S^{n-1}$.
\end{lemma}
\proof Let $x_0\in\partial K$ satisfy that $h_K=\langle u,x_0\rangle$, and hence
$$
h_K(u)-h_K(v)\leq \langle u,x_0\rangle-\langle v,x_0\rangle=
\langle u-v,x_0\rangle\leq\|u-v\|\cdot R.
$$
Since similar argument shows that
$h_K(v)-h_K(u)\leq \|v-u\|\cdot R$, we conclude the lemma.
\proofbox

We also note that if $\mu,\nu$ are Borel probability measures on $S^{n-1}$,
and $f:\,S^{n-1}\to\R$ and $\omega>0$ satisfy that $|f(u)-f(v)|\leq \omega\|u-v\|$ for
$u,v\in S^{n-1}$, then
\begin{equation}
\label{WassersteinLipsitz}
\left|\int_{S^{n-1}}f\,d\mu-\int_{S^{n-1}}f\,d\nu\right|\leq \omega\cdot d_W(\mu,\nu).
\end{equation}

\mbox{ }

\noindent{\bf Proof of Theorem~\ref{VKVCclose} }
Let $d_W(V_K,V_C)=\varepsilon$.
In order to apply Corollary~\ref{cone-volume-diam}, we set
\begin{eqnarray*}
R_0&=& \left\{
\begin{array}{ll}
\left(\frac{n^6}{\delta}\right)^{\frac1{\tau}}&\mbox{ if the action of  $G$ reducible;}\\[1ex]
n& \mbox{ if the action of  $G$ irreducible;}
\end{array} \right.\\[1ex]
r_0&=& \left\{
\begin{array}{ll}
\frac{n^{\frac{n}2}}{5^n}\left(\frac{\delta}{n^6}\right)^{\frac{n-1}{\tau}}
&\mbox{ if the action of  $G$ reducible;}\\[1ex]
\frac1e& \mbox{ if the action of  $G$ irreducible}
\end{array} \right.
\end{eqnarray*}
and  deduce \eqref{r0hKCR0} from  Corollary~\ref{cone-volume-diam}.
In particular, if $u,v\in S^{n-1}$, then
first \eqref{r0hKCR0}, secondly
Lemma~\ref{RKhKLipsitz} and \eqref{r0hKCR0} imply
that if $u,v\in S^{n-1}$, then
\begin{eqnarray*}
|\log h_K(u)-\log h_K(v)|&\leq & \frac{|h_K(u)-h_K(v)|}{r_0}\leq \frac{R_0}{r_0}\cdot\|u-v\|\\
|\log h_C(u)-\log h_C(v)|&\leq & \frac{|h_C(u)-h_C(v)|}{r_0}\leq \frac{R_0}{r_0}\cdot\|u-v\|
\end{eqnarray*}
where
\begin{equation}
\label{R0r0}
\frac{R_0}{r_0}=
\left\{
\begin{array}{ll}
\frac{5^n}{n^{\frac{n}2}}
\left(\frac{n^6}{\delta}\right)^{\frac{n}{\tau}}&\mbox{ if the action of  $G$ reducible;}\\[1ex]
en& \mbox{ if the action of  $G$ irreducible.}
\end{array} \right.
\end{equation}
For
$$
\varepsilon=d_W(V_K,V_C),
$$
we deduce from applying first \eqref{WassersteinLipsitz} and $d_W(V_K,V_C)=\varepsilon$, then
from the Logarithmic Minkowski Inequality Theorem~\ref{logMsym}, and using again
\eqref{WassersteinLipsitz} and $d_W(V_K,V_C)=\varepsilon$
 that
\begin{eqnarray*}
\int_{S^{n-1}}\log h_C\,dV_K&\leq& \int_{S^{n-1}}\log h_C\,dV_C+\frac{R_0}{r_0}\cdot \varepsilon\leq
\int_{S^{n-1}}\log h_K\,dV_C+\frac{R_0}{r_0}\cdot \varepsilon\\
&\leq & \int_{S^{n-1}}\log h_K\,dV_K+\frac{2R_0}{r_0}\cdot \varepsilon.
\end{eqnarray*}
It follows from Theorem~\ref{logMstab}
that  for some $m\geq 1$,
there exist $\theta_1,\ldots,\theta_m>0$ and compact convex sets
$K_1,\ldots,K_m>0$ invariant under $G$ such that
$\sum_{i=1}^m{\rm dim}\,K_i=n$ and
\begin{eqnarray}
\label{c0def}
K_1\oplus\ldots \oplus K_m\subset &K&\subset
\left(1+c_0^n\left(\frac{2R_0}{r_0}\cdot \varepsilon\right)^{\frac1{95n}}\right)(K_1\oplus\ldots \oplus K_m)\\
\nonumber
\theta_1K_1\oplus\ldots \oplus \theta_mK_m\subset &C&\subset
\left(1+c_0^n\left(\frac{2R_0}{r_0}\cdot \varepsilon\right)^{\frac1{95n}}\right)
(\theta_1K_1\oplus\ldots \oplus \theta_mK_m)
\end{eqnarray}
where $c_0>1$ is an absolute constant.

If the action of  $G$ is irreducible, then $m=1$, and hence
$$
\left(1+c_1^n\left(\frac{R_0}{r_0}\cdot \varepsilon\right)^{\frac1{95n}}\right)^{-1}K\subset C
\subset \left(1+c_1^n\left(\frac{R_0}{r_0}\cdot \varepsilon\right)^{\frac1{95n}}\right)K
$$
for some absolute constant $c_1>1$.
In turn, $K,C\subset R_0B^n$ ({\it cf.} \eqref{r0hKCR0}), $R_0=n$
and $\left(\frac{R_0}{r_0}\right)^{\frac1{95n}}<2$ ({\it cf.} \eqref{R0r0})
yield \eqref{KCdilate} as
$$
d_\infty(K,C)\leq R_0\cdot c_1^n\left(\frac{R_0}{r_0}\cdot \varepsilon\right)^{\frac1{95n}}\leq
(2c_1)^n\cdot \varepsilon^{\frac1{95n}}.
$$

Next, let the action of  $G$ be reducible. First, we assume that
\begin{equation}
\label{VKVCclose-epscond}
\varepsilon<c_2^{95n^2}(\delta\tau)^{95n}\left(\frac{\delta}{n^6}\right)^{\frac{96n^2}{\tau}}
\end{equation}
where $c_2\in(0,1)$ is a suitably small absolute constant
 such that if
 $\varepsilon>0$ satisfies \eqref{VKVCclose-epscond}, then
\begin{equation}
\label{delta-tau-eps}
c_0^n\left(\frac{2R_0}{r_0}\cdot \varepsilon\right)^{\frac1{95n}}<\frac{\delta\tau}{4n} \cdot\frac{r_0}{R_0}
\mbox{ \ }(<1)
\end{equation}
holds for the $c_0$ in \eqref{c0def} ({\it cf.} \eqref{R0r0}).
Therefore, on the one hand, we have
$$
\left(1-c_0^n\left(\frac{2R_0}{r_0}\cdot \varepsilon\right)^{\frac1{95n}}\right)
\big((K\cap L_1)\oplus\ldots \oplus (K\cap L_m)\big)\subset K
$$
for $L_i={\rm lin}\,K_i$, $i=1,\ldots,m$,
and, on the other hand, we deduce from \eqref{delta-tau-eps}
and Proposition~\ref{noMinkowskisum} that $m=1$.
In particular,
$$
\left(1-c_3^n\left(\frac{n^6}{\delta}\right)^{\frac{1}{95\tau}}\varepsilon^{\frac1{95n}}\right)K
\subset C\subset
\left(1+c_3^n\left(\frac{n^6}{\delta}\right)^{\frac{1}{95\tau}}\varepsilon^{\frac1{95n}}\right)K
$$
for an suitable absolute constant $c_3>1$, and hence
$K,C\subset R_0B^n$ implies
$$
d_\infty(K,C)\leq R_0\cdot c_3^n\left(\frac{n^6}{\delta}\right)^{\frac{1}{95\tau}}\varepsilon^{\frac1{95n}}.
$$
We conclude Theorem~\ref{VKVCclose} under the condition  \eqref{VKVCclose-epscond}.

Finally, we assume that the condition  \eqref{VKVCclose-epscond} does not hold; namely,
$$
\varepsilon\geq c_2^{95n^2}(\delta\tau)^{95n}\left(\frac{\delta}{n^6}\right)^{\frac{96n^2}{\tau}}.
$$
Since $o\in K,C\subset R_0B^n$, we have
\begin{eqnarray*}
d_\infty(K,C)&\leq &R_0=\left(\frac{n^6}{\delta}\right)^{\frac1{\tau}}\leq
c_2^{-n}(\delta\tau)^{-1}\left(\frac{n^6}{\delta}\right)^{\frac1{\tau}(1+\frac{96n}{95})}\varepsilon^{\frac1{95n}}\\
&\leq &\frac{c_2^{-n}}{\tau}\cdot \delta^{-1}\cdot \left(\frac{n^6}{\delta}\right)^{\frac{2n}{\tau}}
\varepsilon^{\frac1{95n}}\leq \frac{c_2^{-n}}{\tau}\cdot \delta^{\frac{-3n}{\tau}}n^{\frac{12n}{\tau}}
\varepsilon^{\frac1{95n}},
\end{eqnarray*}
proving  Theorem~\ref{VKVCclose}.
\proofbox

\noindent{\bf Proof of Theorem~\ref{VKVCcloseh} }
According to Theorem~\ref{VKsymchar}, there exist convex bodies $K$ and $C$ invariant under $G$
such that $h_1(u)=h_K(u)$ and $h_2(u)=h_C(u)$ for $u\in S^{n-1}$. In turn,
we conclude \eqref{VKVCclosehradii}
from \eqref{r0hKCR0}, and \eqref{VKVCclosehinfty} from
\eqref{r0hKCR0} and \eqref{KCdilate}.\proofbox

After verifying Theorem~\ref{VKVCcloseh}, we consider the case when $\mu_1(S^{n-1})\neq \mu_2(S^{n-1})$.\\

\noindent{\bf Proof of Corollary~\ref{VKVCclosehneq} } We may assume that
$$
1=M=\mu_1(S^{n-1})\leq \mu_2(S^{n-1}).
$$
For $\varepsilon=d_{\rm bL}(\mu_1,\mu_2)\leq 1$, it follows from  \eqref{bLdifference} that
\begin{equation}
\label{mu2eps}
1\leq \mu_2(S^{n-1})\leq 1+\varepsilon.
\end{equation}
We consider the probability measure $\tilde{\mu}_2=\mu_2(S^{n-1})^{-1}\cdot \mu_2$. Since
$d_{bL}(\mu_2,\tilde{\mu}_2)\leq \varepsilon$ by \eqref{bLdifference-lambda}, the triangle inequality yields
$d_{bL}(\mu_1,\tilde{\mu}_2)\leq 2\varepsilon$ by \eqref{mu2eps}, and readily
$$
\tilde{\mu}_2\big(\Psi(L\cap S^{n-1},\delta)\big)\leq (1-\tau)\cdot\frac{{\rm dim}\,L}{n}
$$
for any proper subspace $L$ invariant under $G$. In addition,
$$
\tilde{h}_2=\mu_2(S^{n-1})^{\frac{-1}n}\cdot h_2
$$
is the invariant Alexandrov solution of the Logarithmic Minkowski Problem \eqref{MongeVK-Alexandrov}.

We deduce from Theorem~\ref{VKVCcloseh} that
\begin{eqnarray*}
\|h_1-\tilde{h}_2\|_\infty&\leq & \tilde{\gamma}_0\cdot (2\varepsilon)^{\frac1{95n}}\\[1ex]
r_0\leq &h_1,\tilde{h}_2&\leq \widetilde{R}_0
\end{eqnarray*}
where for some absolute constant $\tilde{c}<1$, we have
\begin{itemize}
\item $\widetilde{R}_0=n$, $r_0=\frac1e$ and  $\tilde{\gamma}_0=\tilde{c}^n\cdot \varepsilon^{\frac1{95n}}$
provided the action of $G$ is irreducible;
\item
$\widetilde{R}_0=\left(\frac{n^6}{\delta}\right)^{\frac1{\tau}}$,
$r_0= \frac{n^{\frac{n}2}}{5^n}\left(\frac{\delta}{n^6}\right)^{\frac{n-1}{\tau}}$
and $\tilde{\gamma}_0=\frac{\tilde{c}^{n}}{\tau}\cdot \delta^{\frac{-3n}{\tau}}n^{\frac{12n}{\tau}}$
provided the action of $G$ is reducible.
\end{itemize}
Therefore,  $h_2=\mu_2(S^{n-1})^{\frac{1}n}\cdot \tilde{h}_2$ and \eqref{mu2eps} imply
Corollary~\ref{VKVCclosehneq} with  $c=2\tilde{c}$ and
$R_0=2\widetilde{R}_0$.
\proofbox

\section{Partial converses Theorem~\ref{KCcloseh} and Theorem~\ref{VKclosemaxh}
of Theorem~\ref{VKVCcloseh}}
\label{secconverse}

In this section, we prove the two partial converses Theorem~\ref{KCcloseh} and Theorem~\ref{VKclosemaxh} of
Theorem~\ref{VKVCcloseh} by verifying Theorem~\ref{KCclose} and Theorem~\ref{VKclosemax}.

Our argument for Theorem~\ref{KCclose} is based on Hug, Schneider \cite{HuS15}, which paper
proved that if $R>0$ and $K$ and $C$ are convex bodies  in $\R^n$  satisfying
$K,C\subset RB^n$,  then
\begin{equation}
\label{HugSchneider}
d_{\rm bL}(S_K,S_C)\leq \tilde{\gamma}(R,n)\cdot \sqrt{d_\infty(K,C)}
\end{equation}
where $\tilde{\gamma}(R,n)>0$ depends on $R$ and $n$.
Theorem~\ref{KCcloseh} directly follows from the following theorem
(see the explanation after \eqref{MongeVK-Alexandrov}).

\begin{theo}
\label{KCclose}
If $R>0$ and $K$ and $C$ are convex bodies  in $\R^n$  satisfying $o\in{\rm int}K,{\rm int}C$ and
$K,C\subset RB^n$,  then
$$
d_{\rm bL}(V_K,V_C)\leq \gamma(R,n)\cdot \sqrt{d_\infty(K,C)}
$$
where $\gamma(R,n)>0$ depends on $R$ and $n$.
\end{theo}
\proof Let $\varepsilon=d_\infty(K,C)\leq R$.
By the symmetry of $K$ and $C$, it is sufficient to prove that
if $f\in {\rm Lip}_1$ with $\|f\|_\infty\leq 1$, then
$$
\int_{S^{n-1}} f\,dV_K-\int_{S^{n-1}} f\,dV_C\leq \gamma(R,n)\cdot \sqrt{\varepsilon}
$$
where $\gamma(R,n)>0$ depends on $R$ and $n$, which is equivalent to say that
\begin{equation}
\label{KCcloseeq}
\int_{S^{n-1}} f\cdot h_K\,dS_K-\int_{S^{n-1}} f\cdot h_C\,dS_C\leq n\gamma(R,n)\cdot \sqrt{\varepsilon}.
\end{equation}
It follows from $d_\infty(K,C)\leq \varepsilon$ that
$$
h_K\leq h_C+\varepsilon.
$$
We deduce from $C\subset R\,B^n$ and Lemma~\ref{RKhKLipsitz} that
$h_C\in {\rm Lip}_R$, and hence  $f\cdot h_C\in {\rm Lip}_{2R}$. For $g=\frac1{2R}\,f\cdot h_C$,
it follows that $g\in {\rm Lip}_{1}$ and $\|g\|_\infty\leq 1$, thus
$\|f\|_\infty\leq 1$, $K\subset R\,B^n$ and
the result \eqref{HugSchneider} by Hug, Schneider \cite{HuS15} yield
\begin{eqnarray*}
\int_{S^{n-1}} f\, h_K\,dS_K-\int_{S^{n-1}} f\, h_C\,dS_C&\leq&
\int_{S^{n-1}} f(h_C+\varepsilon)\,dS_K-\int_{S^{n-1}} f\, h_C\,dS_C\\
&=& \varepsilon\cdot \int_{S^{n-1}} f\,dS_K+\\
&&2R\left(\int_{S^{n-1}} g\,dS_K-\int_{S^{n-1}} g\,dS_C\right)\\
&\leq& \varepsilon\cdot R^{n-1}n\kappa_n+2R\cdot \tilde{\gamma}(R,n)\cdot \sqrt{\varepsilon}.
\end{eqnarray*}
We conclude \eqref{KCcloseeq} from $\varepsilon<2R$, and in turn Theorem~\ref{KCclose}.
\proofbox

Convex bodies whose centroid is the origin and having almost equality in Theorem~\ref{VKsymchar} (ii) were characterized
by B\"or\"oczky, Henk \cite{BoH17}. More precisely, if
$\varepsilon\in(0,\tilde{\varepsilon}_0)$ and the convex body $K\subset \R^n$ have its centroid at the origin,
and satisfies
$$
V_K(L\cap S^{n-1})\geq (1-\varepsilon)\cdot \frac{d}{n}\cdot V(K)
$$
for a linear $d$-space $L$ with $1\leq d<n$, then
\begin{equation}
\label{BoHenk}
 (1-\tilde{\gamma}\cdot\varepsilon^{\frac1{5n}})(C+M)\subset  K\subset C+M
\end{equation}
for some compact convex set $C\subset L^\bot$, and
  complementary  $d$-dimensional compact convex set $M$
where $\tilde{\varepsilon}_0,\tilde{\gamma}>0$ depend on the dimension $n$.

The paper \cite{BoH17} also verified two observations that we need in the sequel. For a convex body $Q$ in $\R^n$,
we write $\sigma(Q)$ to denote the centroid, and $\|x\|_{Q-Q}$ to denote the norm of an $x\in\R^n$ with respect tothe origin symmetric convex body $Q-Q$; namely, $\|x\|_{Q-Q}=\min\{t\geq 0:\,x\in t(Q-Q)\}$.

For convex bodies $K,\widetilde{K}$ in $\R^n$, Lemma~3.4 in \cite{BoH17} says that if
$V(K\Delta \widetilde{K})\leq t\, V(\widetilde{K})$ for $t\in(0,\frac1{4^ne})$, then
\begin{equation}
\label{centroiddistance}
\|\sigma(\widetilde{K})-\sigma(K)\|_{\widetilde{K}-\widetilde{K}}\leq 4nt.
\end{equation}
The second observation, Lemma~3.3 in \cite{BoH17} states that if $z\in\R^n$, then
\begin{equation}
\label{translatediff}
V(\widetilde{K}\Delta (z+\widetilde{K}))\leq 2n
\|z\|_{\widetilde{K}-\widetilde{K}}V(\widetilde{K}).
\end{equation}

The following statement  exhibits why we need a condition of the type of
\eqref{VKVCclosesymcond} in Theorem~\ref{VKVCclose}.

\begin{theo}
\label{VKclosemax}
Let $n\geq 2$, $R>\sqrt{n}$,
 and let
the convex body $K\subset RB^n$ with $V(K)=1$ have its centroid at the origin.
There exist constants
 $\varepsilon_0,\gamma>0$ depending on the dimension $n$,
 such that, if $\varepsilon\in(0,\frac{\varepsilon_0}{R^n})$ and $\delta\in (0,\varepsilon]$ and
 $$
V_K(\Psi(L\cap S^{n-1},\delta))\geq (1-\varepsilon)\cdot \frac{d}{n}
$$
for a linear $d$-space $L$ with $1\leq d<n$, then
 $$
d_\infty(K,C+M)\leq \gamma\,R^{n+1}\varepsilon^{\frac1{5n}}
$$
for some compact convex set $C\subset L^\bot$, and
  complementary  $d$-dimensional compact convex set $M$.

If, in addition, $K$ is invariant under a group
$G\subset O(n)$ leaving $L$ invariant and acting without fixed point on $S^{n-1}$, then
we may assume that $C=K|L^\bot$ and
$M=K|L$.
\end{theo}
\proof We assume that $\varepsilon\in(0,\frac{\varepsilon_0}{R^n})$ where $\varepsilon_0>0$ depending on $n$ is small enough to make the argument work.

We deduce from Lemma~\ref{rRvolume} that $r\,B^n\subset K$ for
$$
r=\frac{n^{\frac{n}2}}{5^nR^{n-1}}.
$$
We plan to cut off a rim from $K$ in order to apply \eqref{BoHenk}. For
$$
\eta=\frac{4\cdot 5^nR^n}{n^{\frac{n}2}}\cdot \delta,
$$
we claim that if $u\in S^{n-1}$ is an exterior normal at $x\in\partial K$
with $x|L\in (1-\eta)(K|L)$, then
\begin{equation}
\label{unotinPsi}
u\not\in \Psi(L\cap S^{n-1},\delta).
\end{equation}
Let $\alpha\in[0,\frac{\pi}2]$, $v\in S^{n-1}\cap L$ and $w\in S^{n-1}\cap L^{\bot}$
such that $u|L=v\cos\alpha$ and $u|L^\bot=w\sin\alpha$, and hence
$u=v\cos\alpha+w\sin\alpha$.

Next let $y\in\partial K$ be such that $v$ is an exterior normal at $y$. Since $x|L\in (1-\eta)(K|L)$, we have
$$
\langle x,v\rangle\leq (1-\eta)\langle y,v\rangle\leq \langle y,v\rangle-\eta r.
$$
It follows that
$$
0\leq h_K(u)-\langle y,u\rangle=\langle x-y,v\cos\alpha+w\sin\alpha\rangle\leq -\eta r\cos\alpha+2R\sin\alpha,
$$
and hence $\tan \alpha\geq \frac{\eta r}{2R}=\frac{\eta\cdot n^{\frac{n}2}}{2\cdot 5^nR^n}$, proving
\eqref{unotinPsi}.

We define
$$
\widetilde{K}=\{x\in K:\,x|L\in (1-\eta)(K|L)\},
$$
thus \eqref{unotinPsi} implies that
\begin{eqnarray}
\nonumber
V_{\widetilde{K}}(L\cap S^{n-1})&\geq &(1-\eta)^n V_K(\Psi(L\cap S^{n-1},\delta))
\geq (1-\gamma_1(R^n\delta+\varepsilon))\cdot \frac{d}{n}\\
\label{VtildeKLlow}
&\geq &
(1-2\gamma_1R^n\varepsilon)\cdot \frac{d}{n}\cdot V(\widetilde{K})
\end{eqnarray}
for $\gamma_1>0$ depending on $n$.
Since $(1-\eta)K\subset\widetilde{K}$, it follows that
$$
V(K\Delta \widetilde{K})\leq \gamma_2R^n V(\widetilde{K})\cdot \varepsilon
$$
for $\gamma_2>0$ depending on $n$.
According to \eqref{centroiddistance} based on \cite{BoH17},
the centroid $\sigma(\widetilde{K})$ of $\widetilde{K}$ satisfies
\begin{equation}
\label{sigmaest}
\|\sigma(\widetilde{K})\|_{\widetilde{K}-\widetilde{K}}\leq 4n\gamma_2R^n\cdot \varepsilon;
\end{equation}
It follows from \eqref{translatediff} based on \cite{BoH17} the convex body $K_0=\widetilde{K}-\sigma(\widetilde{K})$ satifies that
$\sigma(K_0)=o$ and
$$
V(K_0\Delta \widetilde{K})\leq 8n^2\gamma_2R^n V(\widetilde{K})\cdot \varepsilon,
$$
and hence
$$
V_{K_0}(L\cap S^{n-1})\geq V_{\widetilde{K}}(L\cap S^{n-1})-V(K_0\Delta \widetilde{K})\geq
(1-\gamma_3R^n\varepsilon)\cdot \frac{d}{n}\cdot V(K_0)
$$
for $\gamma_3>0$ depending on $n$. We deduce from \eqref{BoHenk} based on \cite{BoH17}
and $V(K_0)\leq 1$ that
there exist some compact convex set $C_0\subset L^\bot$, and
  complementary  $d$-dimensional compact convex set $M_0$ such that
\begin{equation}
\label{BoHenk0}
 (1-\gamma_4R^{\frac1{5}}\cdot\varepsilon^{\frac1{5n}})(C_0+M_0)\subset  K_0\subset C_0+M_0
\end{equation}
where $\gamma_4>0$ depends on the dimension $n$.
Since
$$
K_0+\sigma(\widetilde{K})\subset K\subset (1-\eta)^{-1}(K_0+\sigma(\widetilde{K})),
$$
we deduce from $\widetilde{K}-\widetilde{K}=K_0-K_0$, \eqref{sigmaest} and \eqref{BoHenk0} that
there exist some compact convex set $C\subset L^\bot$, and
  complementary  $d$-dimensional compact convex set $M$ such that
\begin{equation}
\label{BoHenk1}
 (1-\gamma_5R^n\cdot\varepsilon^{\frac1{5n}})(C+M)\subset  K\subset C+M
\end{equation}
where $\gamma_5>0$ depends on the dimension $n$. As $K\subset R\,B^n$, we have
$d_\infty(K,C+M)\leq \gamma_5\,R^{n+1}\varepsilon^{\frac1{5n}}$.

Finally, if $K$ is invariant under a group
$G\subset O(n)$ leaving $L$ (and hence also $L^\bot$) invariant and acting without fixed point on $L\cap S^{n-1}$, then
let $G'\subset O(n)$ be the group whose elements are of the form $\Phi|_L\oplus {\rm id}_{L^\bot}$ for $\Phi\in G$
that acts without non-zero fixed point on $L$,
and  let $G"\subset O(n)$ be the group whose elements are of the form $\Phi|_{L^\bot}\oplus {\rm id}_L$ for $\Phi\in G$
that acts without non-zero fixed point on $L^\bot$.
Now for any $x\in K|L^\bot$, the section
$K\cap(x+L)$ is invariant under  $G'$, and hence the centroid $\sigma(K\cap(x+L))$ of
$K\cap(x+L)$ is invariant under $G'$, which in turn yields that
 $x|L^\bot=\sigma(K\cap(x+L))\in K$. Therefore,
$K|L^\bot=K\cap L^\bot$. Since similar argument implies
$K|L=K\cap L$, we may choose $C=K|L^\bot$ and
$M=K|L$.
\proofbox

\noindent{\bf Proof of Theorem~\ref{VKclosemaxh}: } According to the remarks after \eqref{MongeVK-Alexandrov}, there exists a convex body $K$ invariant under $G$ such that $h=h_K$ and $\mu=V_K$.
Since the centroid of $K$ is invariant under the action of $G$ that does not have non-zero fixed points, it follows that the centroid $\sigma(G)$ of $G$ is the origin.

We deduce from Theorem~\ref{VKclosemax} that there exists $\gamma_0(R,n)>0$ depending on $R$ and $n$ such that
$$
d_\infty(K,C_0+M_0)\leq \gamma_0(R,n)\cdot \varepsilon^{\frac1{5n}}
$$
where $C_0=K|L^\bot$ and
$M_0=K|L$. Rescaling $C_0$ and $M_0$, we obtain convex compact
 $C_1\subset L^\bot$ and $M_1\subset L$ invariant under $G$ such that
$$
d_\infty(K,C_1+M_1)\leq \gamma_1(R,n)\cdot \varepsilon^{\frac1{5n}}
\mbox{ \ and \ }V(C_1+M_1)=1
$$
where $\gamma_1(R,n)>0$ depends on $R$ and $n$. For $Q=C_1+M_1$, it follows from Theorem~\ref{KCclose}
that
$$
d_W(V_K,V_{Q})\leq \gamma(R,n)\cdot \varepsilon^{\frac1{10n}}
$$
where $\gamma(R,n)>0$ depends on $R$ and $n$.

Let $d={\rm dim}\,L$, and let $\varrho>0$ be the maximal radius of a $d$-dimensional ball centered at the origin and contained in $M_1$. We deduce that $d_\infty(\frac{t+\varrho}{\varrho}\,M_1,M_1)\geq t$ for any $t>1$,
and hence
$$
Q_t=\mbox{$\frac{t+\varrho}{\varrho}\cdot M_1+\left(\frac{\varrho}{t+\varrho}\right)^{\frac{d}{n-d}}\cdot C_1$}.
$$
satisfies
\begin{eqnarray*}
d_\infty(K,Q_t)&\geq & t;\\
V(Q_t)&=&1\mbox{ \ and \ }V_{Q_t}=V_{Q};\\
d_W(V_K,V_{Q_t})&\leq &\gamma(R,n)\cdot \varepsilon^{\frac1{10n}}.
\end{eqnarray*}
Therefore, we choose $h_t=h_{Q_t}$ and $\mu_t=V_{Q_t}$.
\proofbox

\noindent{\bf Acknowledgement } We would like to thank Erwin Lutwak, Gaoyong Zhang and Richard Gardner for illuminating discussions.

\end{document}